\documentclass[a4paper,11pt]{article}
\usepackage[margin=0.7in]{geometry}

\usepackage{siunitx}
\usepackage[english]{babel}
\usepackage{amsfonts}
\usepackage{amssymb}
\usepackage{amsmath}
\usepackage{amsthm}
\usepackage{eucal}
\usepackage{mathrsfs}
\allowdisplaybreaks[0]
\usepackage{bm}
\usepackage{graphicx}
\usepackage{caption}
\usepackage{subcaption}
\usepackage{framed}
\usepackage{diagbox}
\usepackage[titletoc,title]{appendix}
\usepackage{tikz}
\usepackage{color}

\usepackage{hyperref}
\usepackage{amsopn}

\usepackage{algorithmicx}
\usepackage[ruled,vlined]{algorithm2e}
\usepackage{algpseudocode,float}
\algblock{ParFor}{EndParFor}

\DeclareMathOperator*{\argmin}{argmin}
\DeclareMathOperator*{\argmax}{argmax}

\author{
Luca Calatroni\thanks{CMAP, \'Ecole Polytechnique, Palaiseau, 91128, Route de Saclay, France
(\href{mailto:luca.calatroni@polytechnique.edu}{luca.calatroni@polytechnique.edu}).} \and Alessandro Lanza\thanks{Department of Mathematics, University of Bologna, Piazza di Porta San Donato 5, Bologna, Italy 
  (\href{mailto:alessandro.lanza2@unibo.it}{alessandro.lanza2@unibo.it},  \href{mailto:monica.pragliola2@unibo.it}{monica.pragliola2@unibo.it}, \href{mailto:fiorella.sgallari@unibo.it}{fiorella.sgallari@unibo.it}).} 
  \and Monica Pragliola\footnotemark[3] 
  \and  Fiorella Sgallari\footnotemark[3]
}

\date{}

\begin{document}
\title{Adaptive parameter selection for weighted-TV \\ image reconstruction problems}

\maketitle

\begin{abstract}
We propose an efficient estimation technique for the automatic selection of locally-adaptive Total Variation regularisation parameters based on an hybrid strategy which combines a local maximum-likelihood approach estimating space-variant image scales with a global discrepancy principle related to noise statistics. We verify the effectiveness of the proposed approach solving some exemplar image reconstruction problems and show its outperformance in comparison to state-of-the-art parameter estimation strategies, the former weighting locally the fit with the data \cite{satv}, the latter relying on a bilevel learning paradigm \cite{Hint1,Hint2}.  
\end{abstract}

\section{Introduction} \label{sec:intro}
In this paper, we are interested in restoring images corrupted by known blur and additive white Gaussian noise (AWGN), i.e. we assume a degradation model of the form $\,g \;{=}\; K u + \epsilon,$
with $g,u,\epsilon \in \mathbb{R}^n$ vectorised forms of the discrete observed image, target image and noise realisation, respectively, and with $K \in \mathbb{R}^{n \times n}$ the blur operator. Such inverse problem is typically ill-posed.

The variational approach to solve ill-posed image restoration problems consists in minimising a composite functional which is the sum of a regularisation term encoding a-priori assumptions on the unknown image $u$ and a data fitting term describing noise statistics.   
A very popular edge-preserving regulariser firstly proposed in \cite{ROF} for image denoising is the Total Variation (TV) semi-norm, while the so-called L$_2$ data term is known to be suited for AWGN. They read
\begin{align}   
\mathrm{TV}(u) &:= 
 \sum_{i=1}^n \|(\mathrm{D}u)_i\|_p \,\;{:=}\;\, 
\sum_{i=1}^n 
\left( (\mathrm{D}_h u)_i^p + (\mathrm{D}_v u)_i^p \right)^{1/p},\quad\;\, p \in\left\{ 1,2\right\} 
\label{TV:def} \tag{TV}\\
\mathrm{L}_2(u) &:=
\frac{1}{2} \, \sum_{i=1}^n (Ku-g)_i^2 , 
\label{L2:def} \tag{L$_2$}
\end{align}

\noindent where $(\mathrm{D}_h u)_i,(\mathrm{D}_v u)_i$ denote the horizontal and vertical discrete gradient components at pixel $i$, respectively. The two instances $p\:{=}1,2$ in \eqref{TV:def} are referred to as \emph{anisotropic} and \emph{isotropic} TV, respectively.
By taking a weighted average of \eqref{TV:def} and \eqref{L2:def}, one gets the TV-L$_2$ model:
\begin{equation}  
\min_{u \:{\in}\: \mathbb{R}^n} 
\left\{ \alpha \text{TV}(u) + \text{L}_2(u) \right\}
\quad\text{or, equivalently,}\quad  
\min_{u \:{\in}\: \mathbb{R}^n} 
\left\{ \text{TV}(u) + \mu \text{L}_2(u) \right\}.
\label{eq:TV-L2}
\tag{TV-L$_2$}
\end{equation}

Both the parameters $\alpha$ and $\mu=1/ \alpha$ are often referred to as \emph{regularisation parameters} since their size weights the amount of the regularisation against the trust in the data. Note that the equivalence of the two formulations in \eqref{eq:TV-L2} allows in fact to use indifferently either of the two models. To estimate an optimal regularisation parameter, several strategies can be used. When the noise level is known, a classical approach is based on the use of the discrepancy principle \cite{discr}, while in blind scenarios, optimisation techniques learning the optimal amount of regularisation from training data can be used, see, e.g., \cite{CalatroniBIL} and the references therein.

In order to overcome the well-known artefacts of TV-based reconstructions, higher-order  \cite{TGV} and/or locally-adaptive anisotropic regularisers \cite{DTV,CMBBE,DTGV} have been proposed in the literature. A simple, though powerful, extension enforcing the TV smoothing to locally adapt to the underlying image structures (such as texture, cartoon\ldots) consists in weighting at any pixel the amount of regularisation \cite{Hint1,Hint2,HintKostas,spatialTV} or data fit \cite{satv}. This reflects in considering two locally-weighted models, referred to as WTV-L$_2$ and TV-WL$_2$, which represent space-variant extensions of the two equivalent formulations of the \eqref{eq:TV-L2} model and read as 
\begin{align}  
&\min_{u\in\mathbb{R}^n} 
\left\{ \mathrm{WTV}(u) + \mathrm{L}_2(u) \right\}, 
&\mathrm{WTV}(u) & {:=} 
\sum_{i=1}^n \alpha_i \| (\mathrm{D} u)_i\|_{p}, 
&\alpha_i&\,{>}\:0 \;\,{\forall}\, i,
\;\, p\in\left\{1,2\right\},
\label{eq:WTV-L2} \tag{WTV-L$_2$}
\\ 
&\min_{u\in\mathbb{R}^n} 
\left\{ \mathrm{TV}(u) + \mathrm{WL}_2(u) \right\}, 
&\mathrm{WL}_2(u) & {:=}
\frac{1}{2} \sum_{i=1}^n \mu_i (Ku-g)_i^2, 
&\mu_i&\,{>}\:0 \;\,{\forall}\, i,
\label{eq:TV-WL2} 
\tag{TV-WL$_2$} 
\end{align}
\noindent where, note, the positive parameters $\alpha_i$ and $\mu_i$ are now space-variant and both control locally the amount of smoothing.
The estimation of parameters $\alpha_i$ in \eqref{eq:WTV-L2} has been done by inferring local geometries \cite{spatialTV} or by means of computationally expensive bilevel-optimisation approaches \cite{Hint1,Hint2}, whereas parameters $\mu_i$ in \eqref{eq:TV-WL2} have been estimated based on the use of a local discrepancy principle \cite{satv}. We remark that despite the aforementioned equivalence of the two approaches in the \eqref{eq:TV-L2} scalar parameter case, the two locally-weighted \eqref{eq:WTV-L2} and \eqref{eq:TV-WL2} models do show significant differences when used for image reconstruction problems as it has been rigorously studied in \cite{HintKostas}. 

\paragraph{Contribution} We propose an image restoration approach  based on an hybrid version of the two space-variant \eqref{eq:WTV-L2} and \eqref{eq:TV-WL2} variational models, with variable regularisation parameters $\alpha_i$ and global fidelity parameter $\mu$ referred to as HWTV-L$_2$ model, see \eqref{eq:HWTV-L2}. We propose a simple yet effective automatic Maximum-Likelihood (ML) estimation procedure of the $\alpha_i$ weights in the WTV regulariser as well as with the use of a standard discrepancy  principle. The statistical prior assumption motivating our ML estimation approach is that image gradients norms are locally drawn from an half-Laplacian distribution with space-variant scale parameters $\alpha_i$. The local closed-form formula obtained by our ML approach is extremely handy and, together with a minimisation algorithm based on the Alternating Directions Method of Multipliers (ADMM), renders our proposal very efficient. The proposed approach outperforms by far both the the classical TV-L$_2$ and the SATV \cite{satv} restoration  methods in terms of standard image quality indexes (ISNR, SSIM) and, compared to recent bilevel optimisation strategies used in \cite{Hint1,Hint2} to estimate local TV weights, it slightly improves the restoration quality while at the same time being much more efficient.

\section{The proposed hybrid space-variant model} 
\label{sec:WTV}

For a given image $g\in\mathbb{R}^n$ corrupted by AWGN and blur generated by a known blur operator $K\in \mathbb{R}^{n\times n}$, we propose the following variational model for image restoration
\begin{equation}  
u^* \;{=}\;\, \argmin_{u\in\mathbb{R}^n} 
\left\{ \mathrm{WTV}(u) + \mu \, \text{L}_2(u) \right\}  \label{eq:HWTV-L2} \tag{HWTV-L$_2$}.
\end{equation}

\noindent We note that in addition to the space-variant parameters $\alpha_i$ contained in the WTV regulariser, a further scalar data weight $\mu>0$ appears in \eqref{eq:HWTV-L2}. This makes our model 
an hybrid version of the two space-variant \eqref{eq:WTV-L2} and \eqref{eq:TV-WL2} models, where \emph{local} parameters $\alpha_i$ describing local image scales in a statistical sense (see Section \ref{sec:ML}) are used together with \emph{global} parameter $\mu$ which codifies the discrepancy w.r.t. the given AWGN level. The redundancy of such parameter is therefore only apparent in \eqref{eq:HWTV-L2} as its value is computed depending on the global noise statistics in comparison with the local regularisation strength encoded by the parameters $\alpha_i$.

\subsection{Statistical derivation}

We now justify the choice of the space-variant WTV regulariser in \eqref{eq:WTV-L2} by means of statistical arguments.

A common paradigm in image restoration is the Maximum A Posteriori (MAP) approach by which the restored image is obtained as a global minimiser of the negative log-likelihood distribution of the observed image $g$ given the blurring operator $K$ combined with some prior PDF $P(u)$ on the unknown target image $u$.
In formulas:
\begin{equation}
u^* \in \argmax_{u \in \mathbb{R}^n}\;
P(u|g;K) 
\;{=}\;
\argmin_{u \in \mathbb{R}^n} \;
\left\{ \,-\log P(g|u;K)-\log P(u)
\, \right\},
\label{eq:map2}
\end{equation}
where equality comes from the Bayes' formula after dropping the normalisation term $P(g)$.
The terms $P(g|u;K)$ and $P(u)$ in \eqref{eq:map2} are commonly referred to as the \emph{likelihood} and \emph{prior} distribution, respectively: they encode available information on the statistics of the noise and on the solution we seek, respectively.

A standard prior model for the gradient magnitude of the unknown image $u$ is TV-Gibbs prior so that 
$
P(u)= \frac{1}{Z} \exp \left(-\alpha \sum_{i=1}^n \| (\mathrm{D}u)_i \|_p\right)
$, where $p\in\{1,2\}$.
As pointed out in \cite{Vip},
%
%
such choice can be equivalently interpreted by saying that each $\|(\mathrm{D} u)_i\|_{p}$ distributes according to
an half-Laplacian PDF with parameter $\alpha>0$.
The use of a one-parameter distribution may be restrictive in  modelling images with local properties at different scales (edges, texture\ldots). To allow more flexibility, in \cite{Vip} a space-variant model where gradient norms distribute according to a half-Laplacian distribution with \emph{locally varying} scale parameter $\alpha_i>0$ has been proposed. The prior associated to such choice is
\begin{equation}  \label{prior:WTV}
P(u)= \frac{1}{Z} \,\text{exp}\,\left(-\,\sum_{i = 1}^{n}\alpha_i || (\mathrm{D} u)_{i}||_{p}\right) = \frac{1}{Z} \,\text{exp}\,\left(- \mathrm{WTV}(u)\right),\quad p\in\left\{1,2\right\},
\end{equation}
where WTV is the regulariser defined in \eqref{eq:WTV-L2}.

For the sake of completeness, we recall that in the case of AWGN the likelihood term in \eqref{eq:map2} takes the following form
\begin{equation}
P(g|u;K)=\prod_{i=1}^{n}\,\frac{1}{\sqrt{2\pi}\sigma}\,\text{exp}\bigg(-\frac{(Ku-g)_{i}^{2}}{2\sigma^{2}}\,\bigg)=\frac{1}{W}\,\text{exp}\bigg(-\frac{\lVert Ku-g\rVert_{2}^{2}}{2\sigma^{2}}\,\bigg),
\label{eq:g_like}
\end{equation}
where $\sigma>0$ denotes the AWGN standard deviation and $W > 0$ is a normalisation constant. By plugging the expression \eqref{prior:WTV} for $P(u)$ and \eqref{eq:g_like} $P(g|u;K)$ in \eqref{eq:map2}, we derive the variational model \eqref{eq:WTV-L2}.

In our modelling, in order to describe local image features together with global noise discrepancy, we further weight the data fitting term by a global parameter $\mu$, thus obtaining the hybrid reconstruction model \eqref{eq:HWTV-L2}.

\section{ADMM optimisation \& automatic parameter selection}


In order to solve numerically the image restoration problem \eqref{eq:HWTV-L2}, we use in the following an ADMM-based algorithm combined with an adaptive estimation procedure of model parameters along the iterations. To do so, we introduce two auxiliary variables $w \in \mathbb{R}^n$ and $t \in \mathbb{R}^{2n}$ and rewrite the model in the following constrained form:
\begin{eqnarray}\label{eq:PM_ADMM_a}
\{ \, u^*,w^*,t^* \}
&\:\;{\leftarrow}\;\:&
\argmin_{u,w,t}
\bigg\{ \:
\sum_{i = 1}^{n}\alpha_i \left\| t_i \right\|_{p}
\;{+}\;
\frac{\mu}{2} \, \| w \|_2^2
\: \bigg\} \\
&\text{subject to}&\quad w \;{=}\; K u - g, \;\:
t \;{=}\; \mathrm{D}u. \notag
\end{eqnarray}
We first define the augmented Lagrangian functional:
\begin{eqnarray}\nonumber
\mathcal{L}(u,w,t;\rho_w,\rho_t;\alpha_1,\ldots,\alpha_n, \mu)&:=& \sum_{i = 1}^{n} \alpha_i \left\|  t_i \right\|_{p}+\frac{\mu}{2} \| w \|_2^2 - \rho_t^{T}(t - \mathrm{D}u) +\frac{\beta_t}{2} \| t - \mathrm{D} u \|_2^2\\
\label{eq:PM_AL}
&-&\rho_w^{T}(w - (Ku-g))+ \frac{\beta_w}{2}\| w - (Ku-g) \|_2^2, 
\end{eqnarray}
\noindent where $\beta_w, \beta_t > 0$ are scalar penalty parameters and $\rho_w \in \mathbb{R}^n$, $\rho_t \in \mathbb{R}^{2n}$
are the vectors of Lagrange multipliers. 
The solution $(u^*,w^*,t^*)$ of problem (\ref{eq:PM_ADMM_a}) is a saddle point for $\mathcal{L}$ in (\ref{eq:PM_AL}). Hence, we can alternate a minimisation step with respect to the primal variables $t,u,w$ with a maximisation step with respect to the dual variables $\rho_t,\rho_w$, in combination with an iterative update of the space variant parameters $\alpha_i$ and $\mu$, which hence will be denoted by $\alpha^{(k)}_i$ and $\mu^{(k)}$. In particular, for what concerns $\alpha_i^{(k)}$ we use the easy ML estimation strategy described next, whereas for $\mu^{(k)}$ we will rely on a global discrepancy principle 

%
\paragraph{Primal variables update.} The three primal sub-problems can be solved efficiently and in closed-form by simple shrinkage/projection operators (in both cases $p=1,2$) and linear system solvers - see \cite{TVp,Vip}. More in details, the sub-problem with respect to the $t$ primal variable, after some algebraic manipulations, can be written as,
\begin{eqnarray}\label{eq:PM_ADMM_a}
t^{(k+1)}&\leftarrow&\argmin_{t}\sum_{i = 1}^{n}\bigg\{ \:
\alpha_i^{(k)} \left\| t_i \right\|_{p} \;{+}\; \frac{\beta_t}{2} \bigg\| t_i - \bigg((\mathrm{D}u^{(k)})_i + \frac{1}{\beta_t}(\rho_t^{(k)})_i\bigg)  \bigg\|_2^2\bigg\}.
\end{eqnarray}
Denoting by,
\begin{equation}
q_i^{(k)} = (\mathrm{D}u^{(k)})_i + \frac{1}{\beta_t}(\rho_t^{(k)})_i\,\in \mathbb{R}^2,
\end{equation}
the solution of each one-dimensional separable problem is given by
\begin{equation}
t_i^{(k+1)} = q_i \max\left(1 - \frac{\alpha_i^{(k)}}{\beta_t\lVert q_i^{(k)} \rVert_p} ,0\right)\,,\quad p \in\{1,2\}\,,\quad i=1,\ldots,n.
\end{equation}

Introducing
\begin{equation}
z^{(k)} \;{=}\;\: Ku^{(k)} -\, g \: + \, \frac{1}{\beta_w} \, \rho_w^{(k)} \; ,
\label{eq:v_def}
\end{equation}
we have that the updating formula for $w$ reads:
\begin{equation}
\label{eq:sub_r_q2_q2_sol}
w^{(k+1)}= \frac{\beta_w}{\mu^{(k)} + \beta_w} z^{(k)}.
\end{equation}
Imposing a first order optimality condition with respect to the primal variable $u$, leads to the following linear system,
\begin{equation}
\left(
\mathrm{D}^T \mathrm{D}
+ \frac{\beta_w}{\beta_t} K^T K
\right)
u
=
\mathrm{D}^T \left( t^{(k+1)} - \frac{1}{\beta_t} \rho^{(k)}_t \right)
+
\frac{\beta_w}{\beta_t} K^T \left( w^{(k+1)} - \frac{1}{\beta_w} \rho^{(k)}_w + g  \right)
\: ,
\label{eq:sub_u_sol}
\end{equation}
that can be solved since the coefficient matrix has full rank - see, e.g, \cite{DTV}.
%


\paragraph{Parameters update.} \label{sec:ML}
In Algorithm \ref{alg:ADMM} both the local space-variant parameters $\alpha_i$ and the global parameter $\mu$ are updated along the iterations. This is a standard strategy for this type of optimisation problems (see, e.g., \cite{DiscrepancyADMM,DTV}), especially in the case of a cheap update of parameters adapting to the image quality improvement as the one considered in the following. Despite the iterative change in the expression of the cost functional corresponding to such update, we remark that nonetheless we observed empirical convergence of the algorithm. A theoretical proof of such result is left for future research. 

For any pixel $i=1,\ldots, n$, we consider the set $\mathcal{S}_i:=\{x_{i,j}\}_{j=1}^N$ with $x_{i,j}=\|(\mathrm{D} u^{(k)})_j\|_{p}$, where $(\mathrm{D} u^{(k)})_j$ are gradients in the square neighbourhood $\mathcal{N}_{i}^{r}$ centred in $i$ with side $2r+1$ whose norm is drawn from a half-Laplacian distribution with scale parameter $\alpha_i$. The likelihood function of $\alpha_i$ thus reads:
\begin{equation}
\mathcal{L}(\alpha_i;\mathcal{S}_i) \;=\;\prod_{x_{i,j}\in\mathcal{S}_i} P(x_{i,j};\alpha_i)=\prod_{j=1}^{N} P(x_{i,j};\alpha_i)\;=\;\alpha_i^{N}\exp\bigg(-\sum_{j=1}^{N} \alpha_i x_{i,j}\bigg).
\label{likelihood}
\end{equation}
We now look for an $\alpha_i>0$ maximising $\mathcal{L}$, or equivalently, minimising $\mathcal{F}(\alpha_i;\mathcal{S}_i) :=-\log~\mathcal{L}(\alpha_i;\mathcal{S}_i)$.
By imposing first order optimality on $\mathcal{F}$ with respect to $\alpha_i$, we obtain the closed formula:
\begin{equation}
\alpha_i = \bigg(\frac{1}{N}\sum_{j=1}^{N}x_{i,j}\bigg)^{-1},
\label{est}
\end{equation}
\noindent which can be handily updated along the iterations $k\geq 0$ to estimate the local regularisation parameters $\alpha^{(k)}_i$ at each pixel $i=1,\ldots,n$ by taking as samples $x^{(k)}_{i,j}=\|(\mathrm{D} u^{(k)})_j\|_{p}$, $j=1,\ldots,N$ i.e. the norms of the image gradients in the neighbourhood $\mathcal{N}_i^{r}$.
We remark that the estimates $\alpha_i^{(k)}$ in (\ref{est}) can be  efficiently computed based on 2D convolution (realised by a fast 2D discrete transform) of the map of gradient norms with a square $(2r+1) \times (2r+1)$ averaging kernel.

The  parameter $\mu$ is updated along the iterations so as to fulfil the global discrepancy principle as described in \cite{DiscrepancyADMM}: we ask each iterate $u^{(k)}$ to satisfy the condition
$|| K u^{(k)} - g ||_2  \leq \delta :=\tau \sigma \sqrt{n}$, where $\sigma$ is the noise standard deviation and the parameter $\tau \approx 1$ is set a priori (see Section \ref{sec:ex} for some experiments describing the sensitivity of the model to this parameter).
In particular, recalling the definition of $z^{(k)}$ given in \eqref{eq:v_def}, the update reads:
\begin{align}
\|  z^{(k)} \|_2 \leq \delta \quad
&\!\!\Longrightarrow\!\! \quad
\mu^{(k+1)} = 0, \label{eq:update_mu} \\
\| z^{(k)} \|_2 > \delta \quad
&\!\!\Longrightarrow\!\!  \quad
\mu^{(k+1)} = \beta_w \big( \| \, z^{(k)} \|_2 / \delta - 1 \big).  \notag
\end{align}

The ADMM pseudo-code is reported in Algorithm \ref{alg:ADMM}.
%

\begin{minipage}[t]{0.85\textwidth}\null
\begin{algorithm}[H]
\caption{ADMM for HWTV-L$_2$ model}
\label{alg:ADMM}
\SetKwInOut{Input}{Input}
\SetKwInOut{Parameters}{Parameters}
\SetKwInOut{Identification}{Identification}
\SetKwInOut{Initialization}{Initialization}
\SetKwInOut{Update}{Update}
\SetKw{Init}{Initialise}
\SetKw{Est}{Estimate}
\Input{observed image $g\in\mathbb{R}^n$;}
\Parameters{$r>0,~ \tau\approx 1,~ \beta_t,~\beta_w >0$; }
\Init{$u^{(0)}=g$, $\rho_w^{(0)}=\rho_t^{(0)}=0$;}\\
\For{k=0,1,\ldots\quad \emph{until convergence}}
{\textbf{update parameters:}\\

$$
\begin{aligned}
& \alpha_i^{(k)} \text{ by } \eqref{est} \text{ for every }i=1,\ldots,n,  \\
& \mu^{(k)} \text{ by } \eqref{eq:update_mu},\\
\end{aligned}
$$

\textbf{update primal variables:}\\

$$
\begin{aligned}
&t^{(k+1)}=\argmin_t \mathcal{L}(u^{(k)},w^{(k)},t;\rho_w^{(k)},\rho_t^{(k)};\alpha^{(k)}_1,\ldots,\alpha^{(k)}_n,\mu^{(k)})\\
&w^{(k+1)}=\argmin_r \mathcal{L}(u^{(k)},w,t^{(k+1)};\rho_w^{(k)},\rho_t^{(k)};\alpha^{(k)}_1,\ldots,\alpha^{(k)}_n,\mu^{(k)})\\
&u^{(k+1)}=\nolinebreak[4]
\argmin_u\mathcal{L}(u,w^{(k+1)},t^{(k+1)};\rho_w^{(k)},\rho_t^{(k)};\alpha^{(k)}_1,\ldots,\alpha^{(k)}_n,\mu^{(k)})\\
\end{aligned}
$$

\textbf{update dual variables:}\\

$$
\begin{aligned}
&\rho_w^{(k+1)}=\rho_w^{(k)}-\beta_w \big(w^{(k+1)}-(K u^{(k+1)}-g)\big),\\
&\rho_t^{(k+1)}=\rho_t^{(k)}-\beta_t \big(t^{(k+1)}-\mathrm{D} u^{(k+1)}\big),\\
\end{aligned}
$$

}
\Return{$u^{*}=u^{(k+1)}$.}
\end{algorithm}
\end{minipage}


\section{Numerical results}
\label{sec:ex}

In this section we report some numerical results obtained by solving the image reconstruction model \eqref{eq:HWTV-L2} via the ADMM Algorithm \ref{alg:ADMM} with fixed penalty parameters $\beta_t=20$ and $\beta_w=100$. In our experiments we observed that the convergence properties of the algorithm are not affected by this choice, if not in terms of convergence speed. The value $r>0$ denotes the radius of the neighbourhoods $\mathcal{N}_i^r$ defined in Sect. \ref{sec:ML} and used to estimate the space-variant parameters $\alpha_i$. Denoting by $u\in [0,1]^n$ the ground-truth image, we assess the quality of the reconstruction $u^*$ by means of 
the Improved Signal-to-Noise Ratio $\mathrm{ISNR}(g,u,u^*) := 10\log_{10} (\|g-u\|_2^2 / \|u^*-u\|_2^2 )$
and in terms of the Structural Similarity Index (SSIM). We compare our results with the ones obtained by the standard \eqref{eq:TV-L2} model, the SATV approach \cite{satv} based on coupling the \eqref{eq:TV-WL2} model with a local discrepancy-based procedure for automatically selecting the parameters $\mu_i$ \footnote{We used the MATLAB code available at: \url{https://www.math.hu-berlin.de/~hp_hint/software/satv.html}.} and the bilevel learning strategy used in \cite{Hint1,Hint2} to estimate the parameters $\alpha_i$ of \eqref{eq:WTV-L2} model 
via a nested optimisation procedure.  

\begin{figure}[!t]
	\centering
	\begin{subfigure}{0.22\textwidth}
	\centering
		\includegraphics[scale = 0.35]{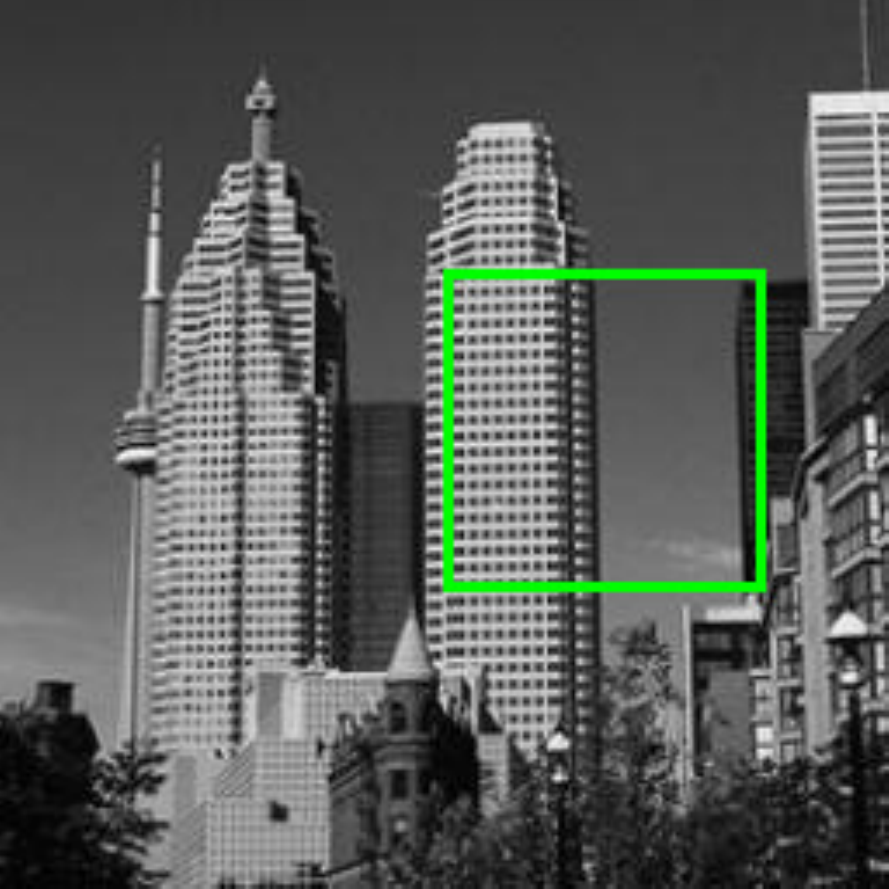}
	  \caption{Original $u$.}
		\label{fig:skyscraper_or}
	\end{subfigure}
	\begin{subfigure}{0.22\textwidth}
	\centering
		\includegraphics[scale = 0.35]{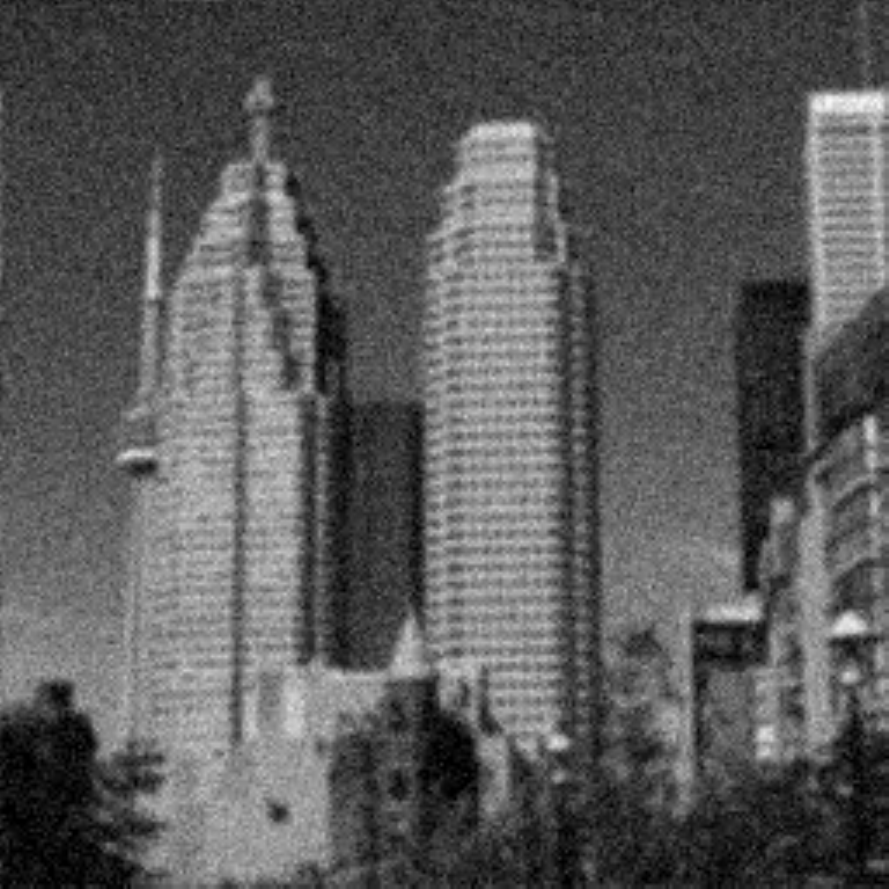}
	  \caption{Corrupted $g$.}
		\label{fig:skyscraper_cor}
	\end{subfigure}
	\caption{Ground truth (with detail) and noisy version corrupted with AWGN with $\sigma=0.05$.}
 \end{figure}

\begin{figure}[!t]
	\centering
	\begin{subfigure}{0.24\textwidth}
	\centering
		\includegraphics[scale = 0.34]{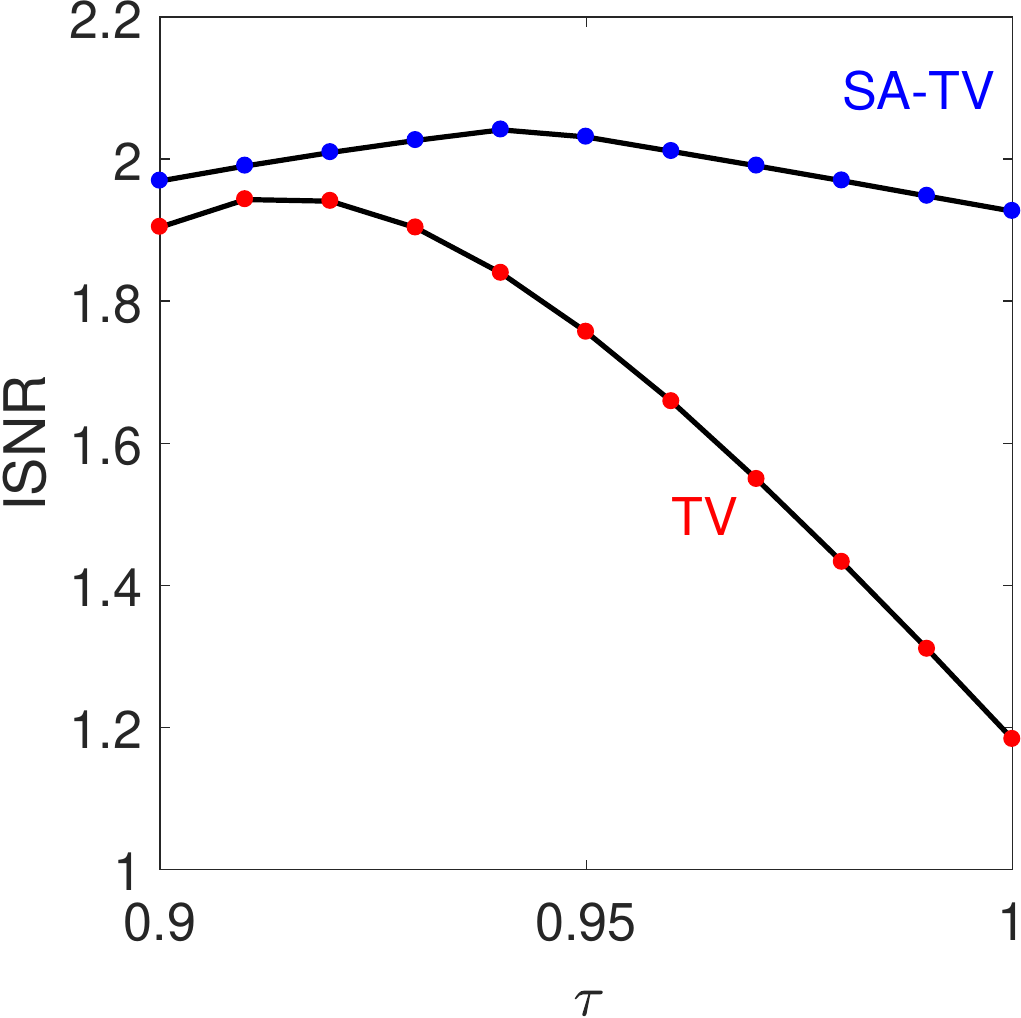}
		\caption{ISNR vs $\tau$}
		\label{fig:ISNR_TVSATV}
	\end{subfigure}
	\begin{subfigure}{0.24\textwidth}
	\centering
	   \includegraphics[scale = 0.34]{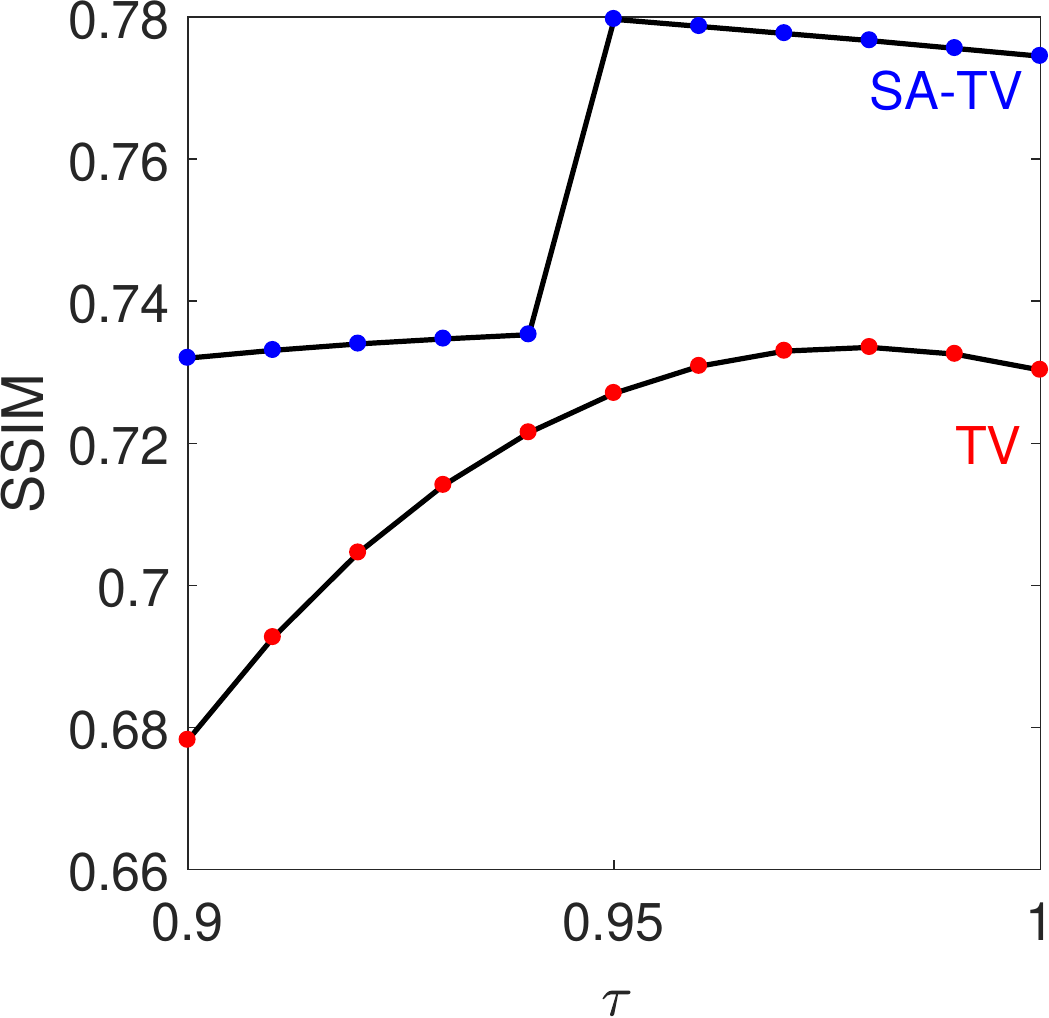}
		\caption{SSIM vs $\tau$}
		\label{fig:SSIM_TVSATV}
	\end{subfigure}
    \begin{subfigure}{0.24\textwidth}
    \centering
	   \includegraphics[scale = 0.36]{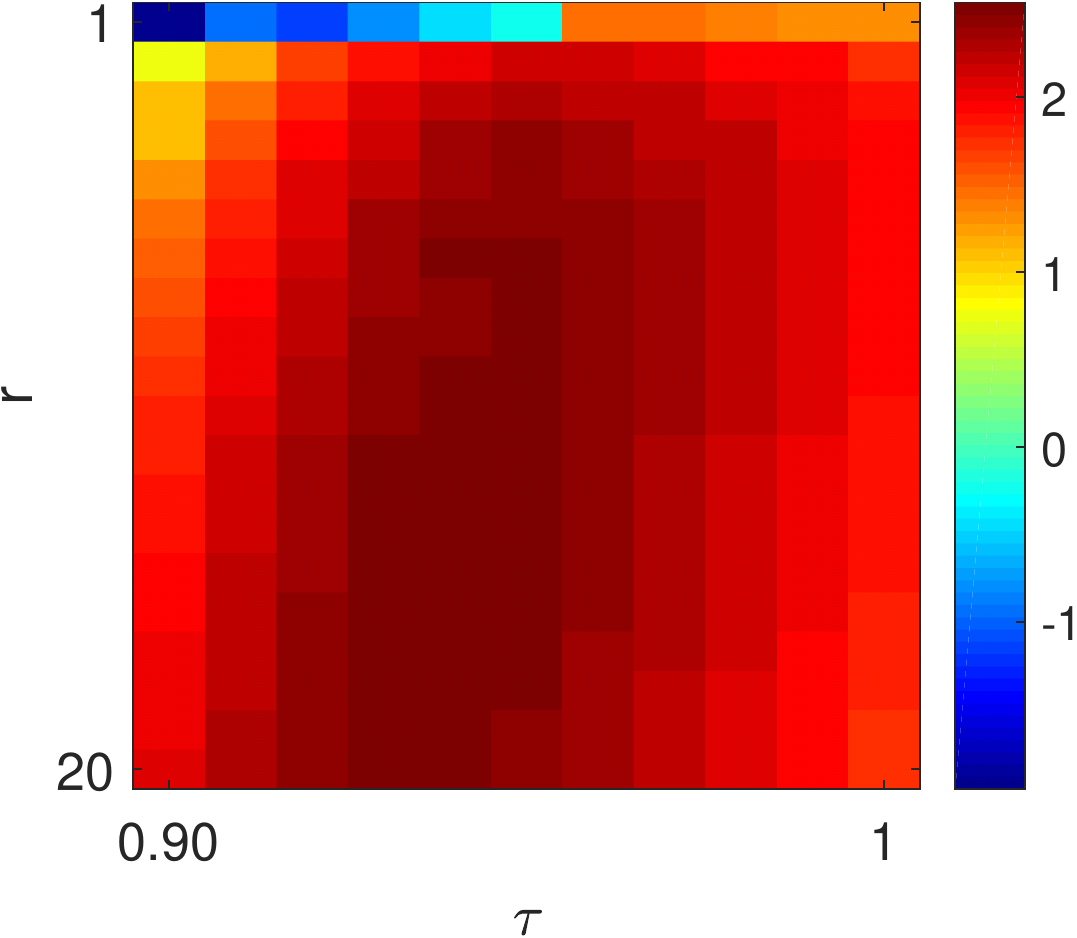}
		\caption{ISNR$(\tau,r)$}
		\label{fig:ISNR_SV}
	\end{subfigure}
    \begin{subfigure}{0.24\textwidth}
    \centering
	   \includegraphics[scale = 0.36]{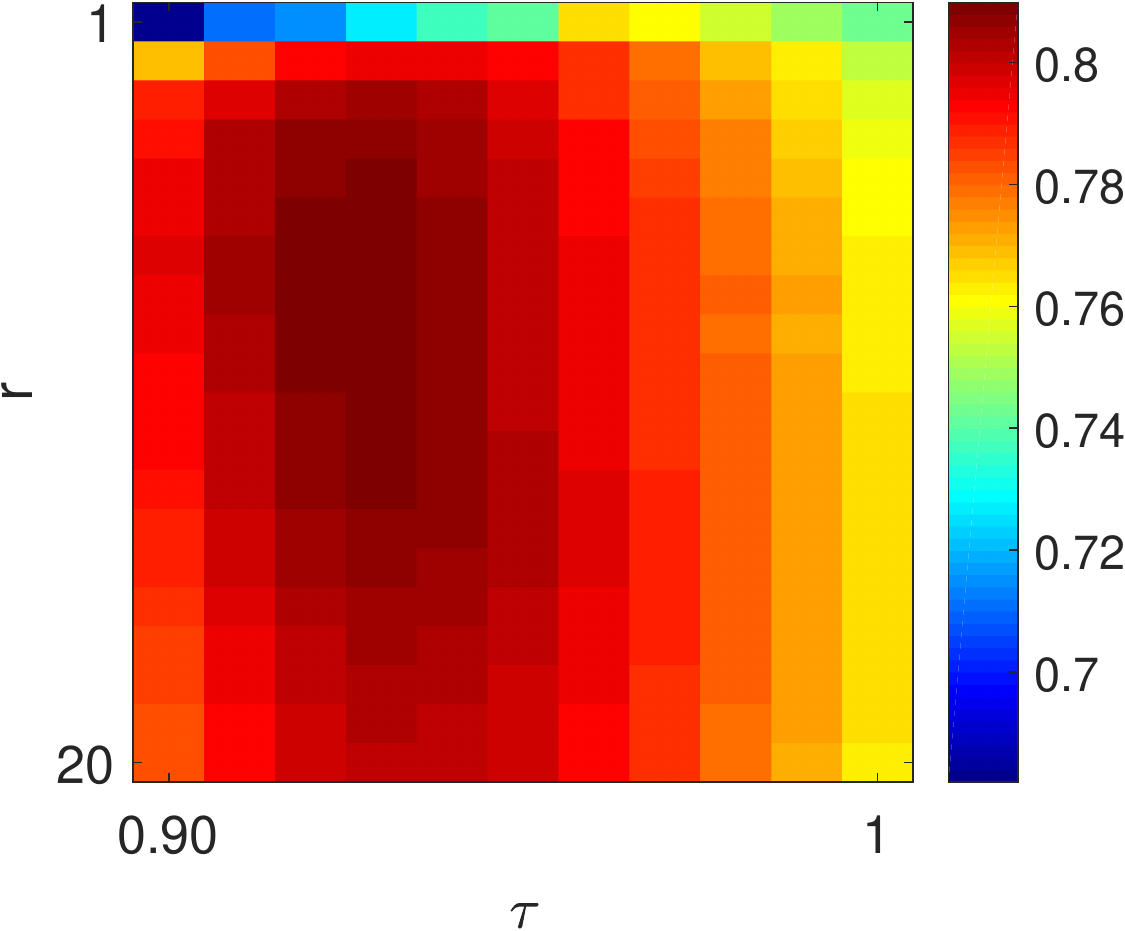}
		\caption{SSIM$(\tau,r)$}
		\label{fig:SSIM_SV}
	\end{subfigure}
   \caption{ISNR \ref{fig:ISNR_TVSATV} and SSIM \ref{fig:SSIM_TVSATV} values reached for different values of $\tau$ by applying TV-L$_2$ and SATV to the restoration of \texttt{skyscraper} test image in Fig.\ref{fig:skyscraper_or}(bottom). For the same image, ISNR  \ref{fig:ISNR_SV} and SSIM \ref{fig:SSIM_SV} values achieved by HWTV-L$_2$ method for different values of $\tau$ and $r$.}
    	\label{fig:ssimisnr}
\end{figure}

\begin{figure}[!t]
	\centering
		\begin{subfigure}{0.24\textwidth}
		\centering
		\includegraphics[scale = 0.34]{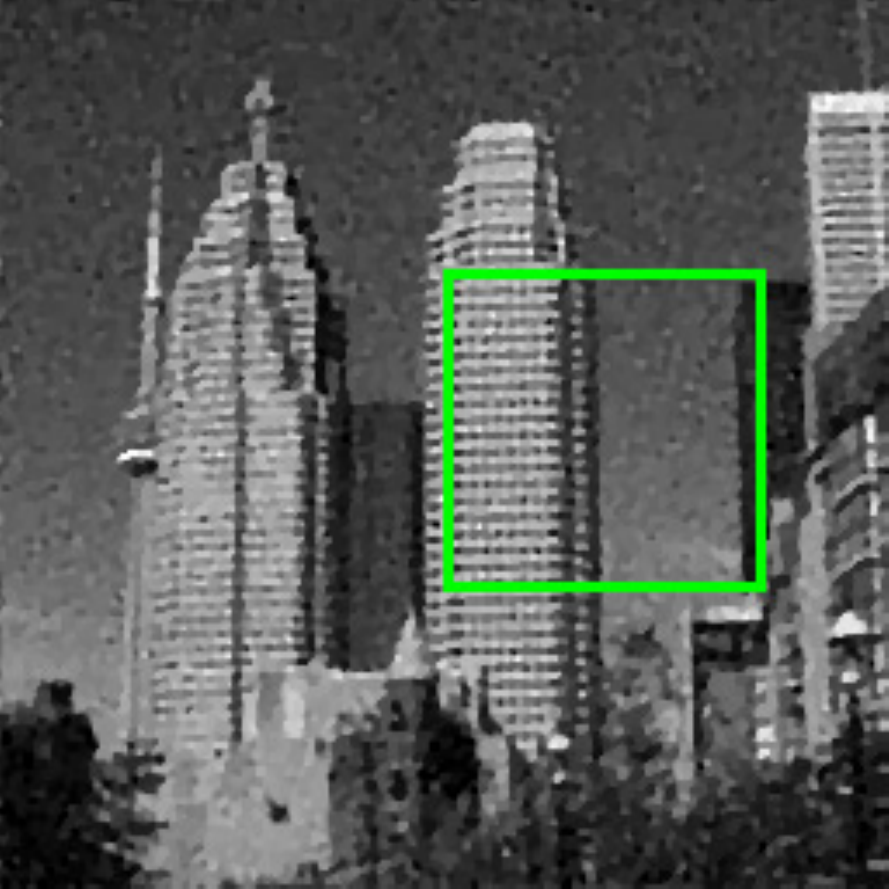}
		\caption{TV-L$_2$}
		\label{fig:TVsky}
		\end{subfigure}
		\begin{subfigure}{0.24\textwidth}
		\centering
		\includegraphics[scale=0.34]{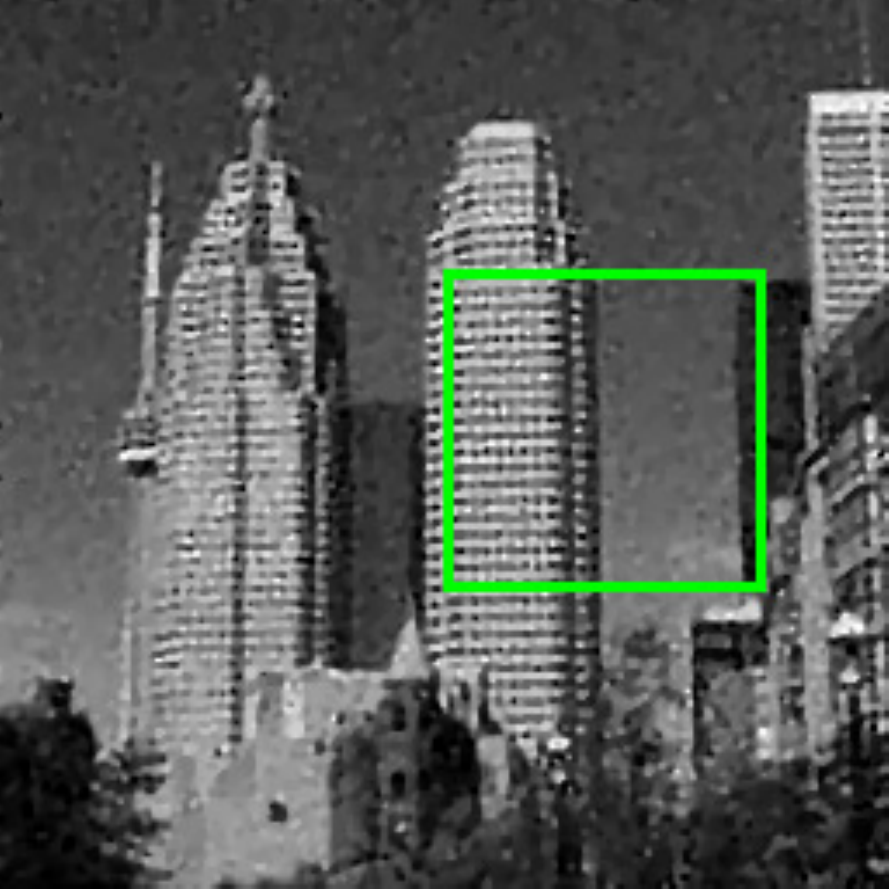}
		\caption{SATV}
		\label{fig:SATVsky}
		\end{subfigure}
		\begin{subfigure}{0.24\textwidth}
		\centering
		\includegraphics[scale = 0.34]{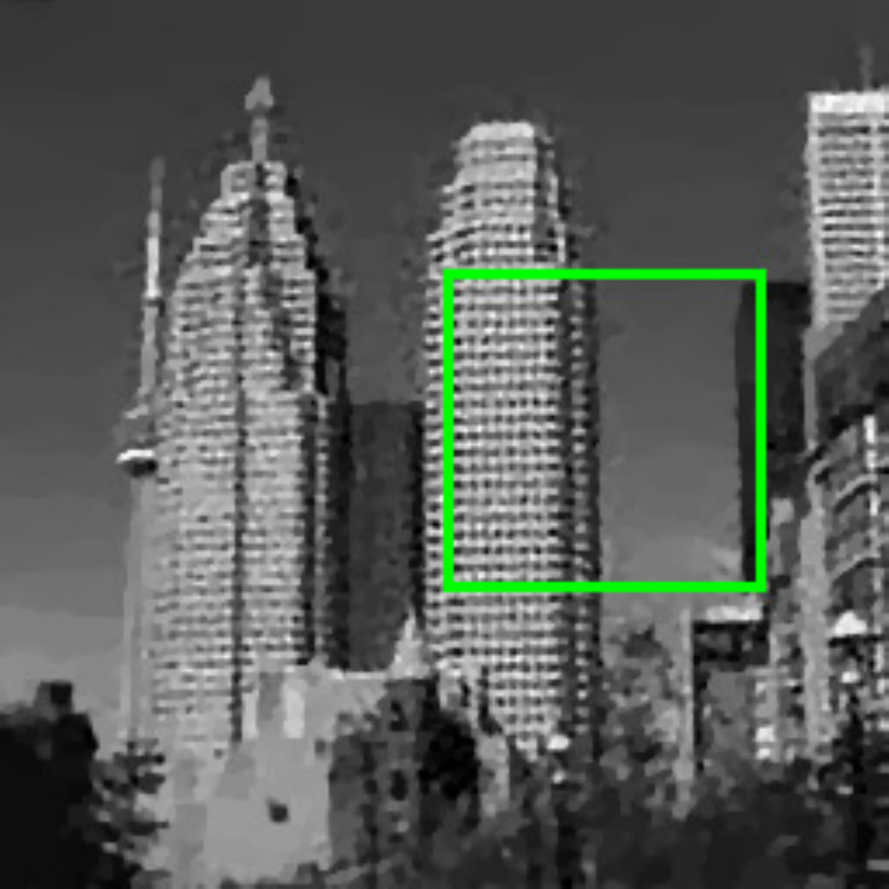}
		\caption{HWTV-L$_2$}
		\label{fig:WTVsky}
		\end{subfigure} 
		\begin{subfigure}{0.24\textwidth}
		\centering
		\includegraphics[scale = 0.20]{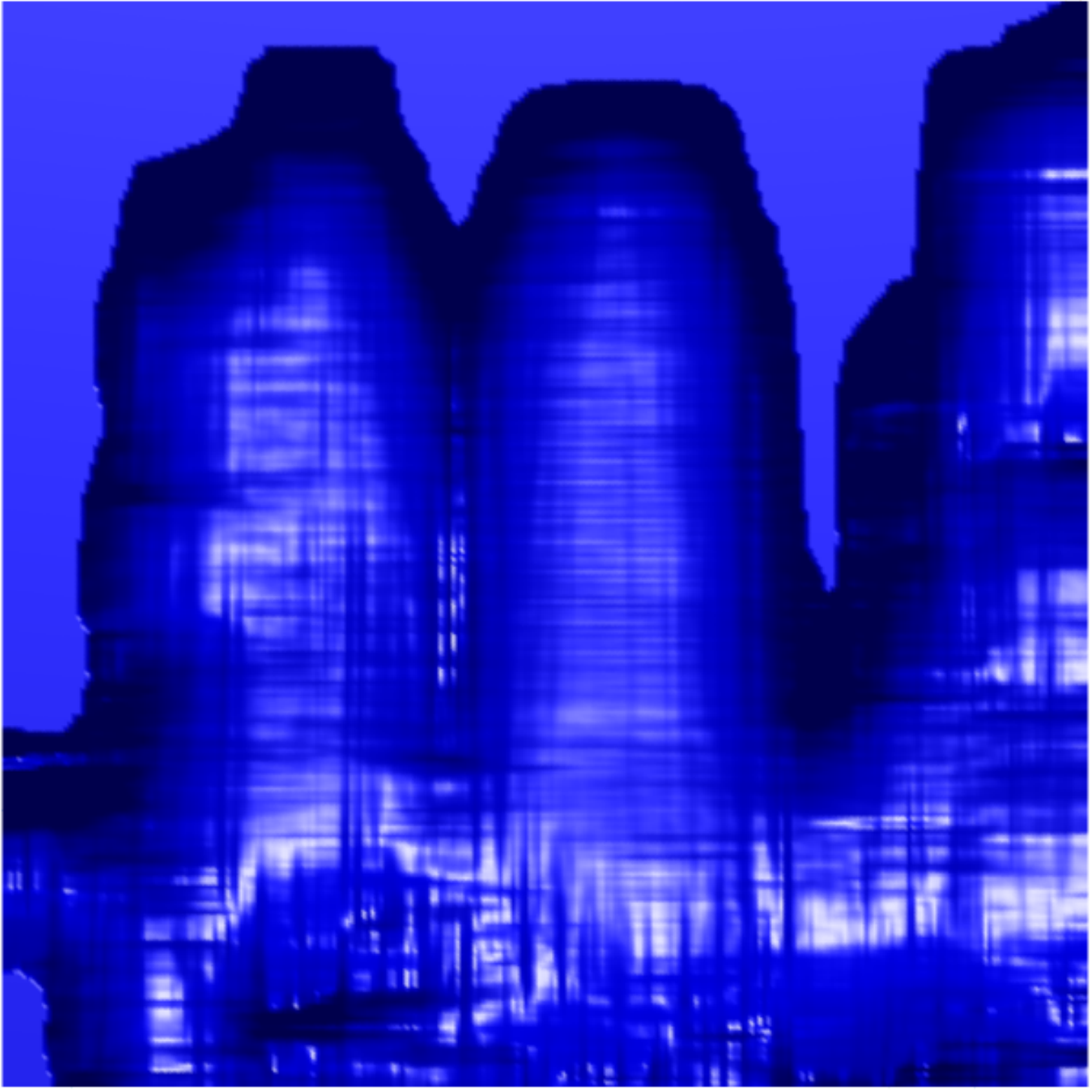}
		\caption{$\alpha$}
        \label{fig:alpha}
		\end{subfigure}\\
		\begin{subfigure}{0.24\textwidth}
		\centering
		\includegraphics[scale = 0.25]{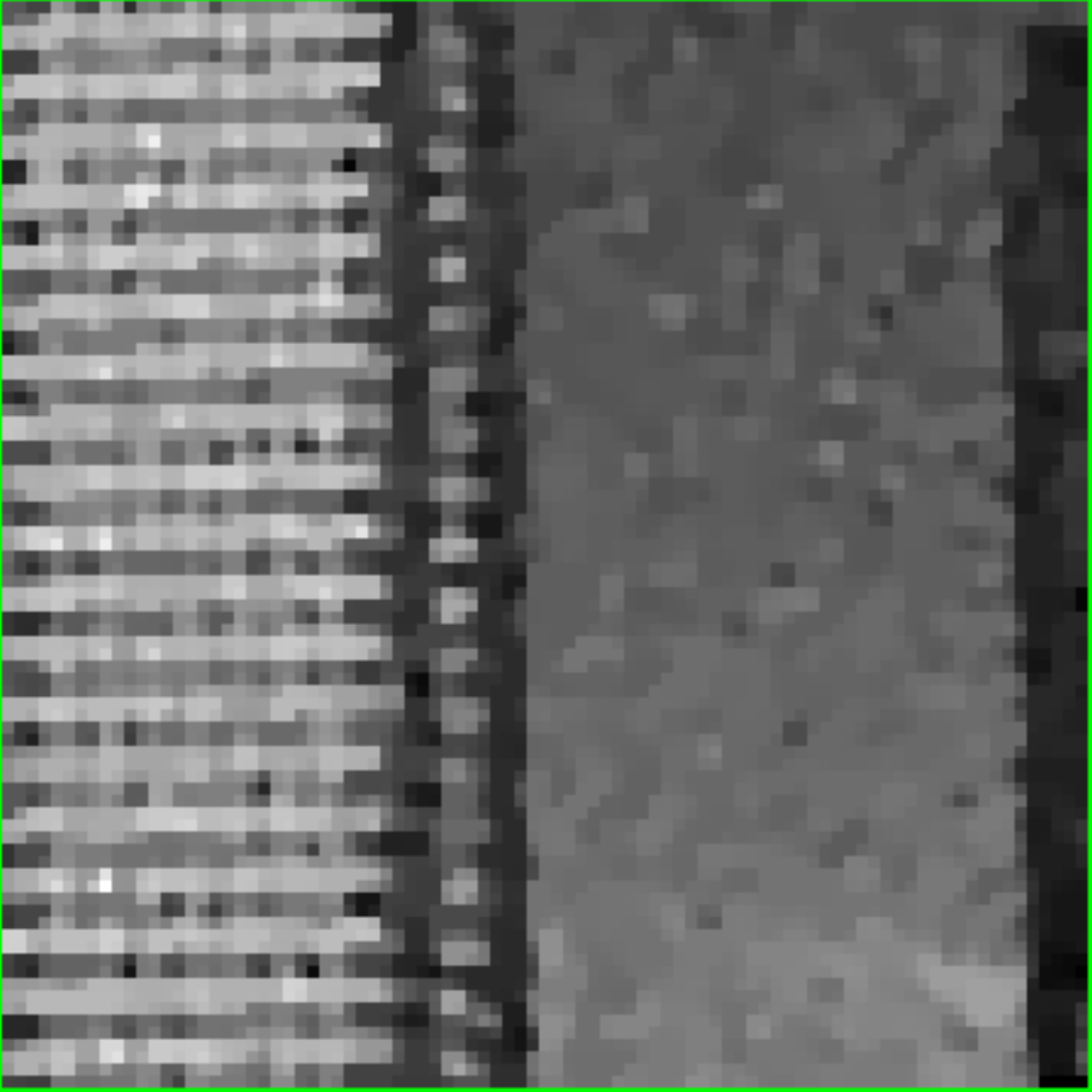}
		\caption{TV-L$_2$ (zoom)}
		\label{fig:TVsky_zoom}
		\end{subfigure}
		\begin{subfigure}{0.24\textwidth}
		\centering
				\includegraphics[scale=0.25]{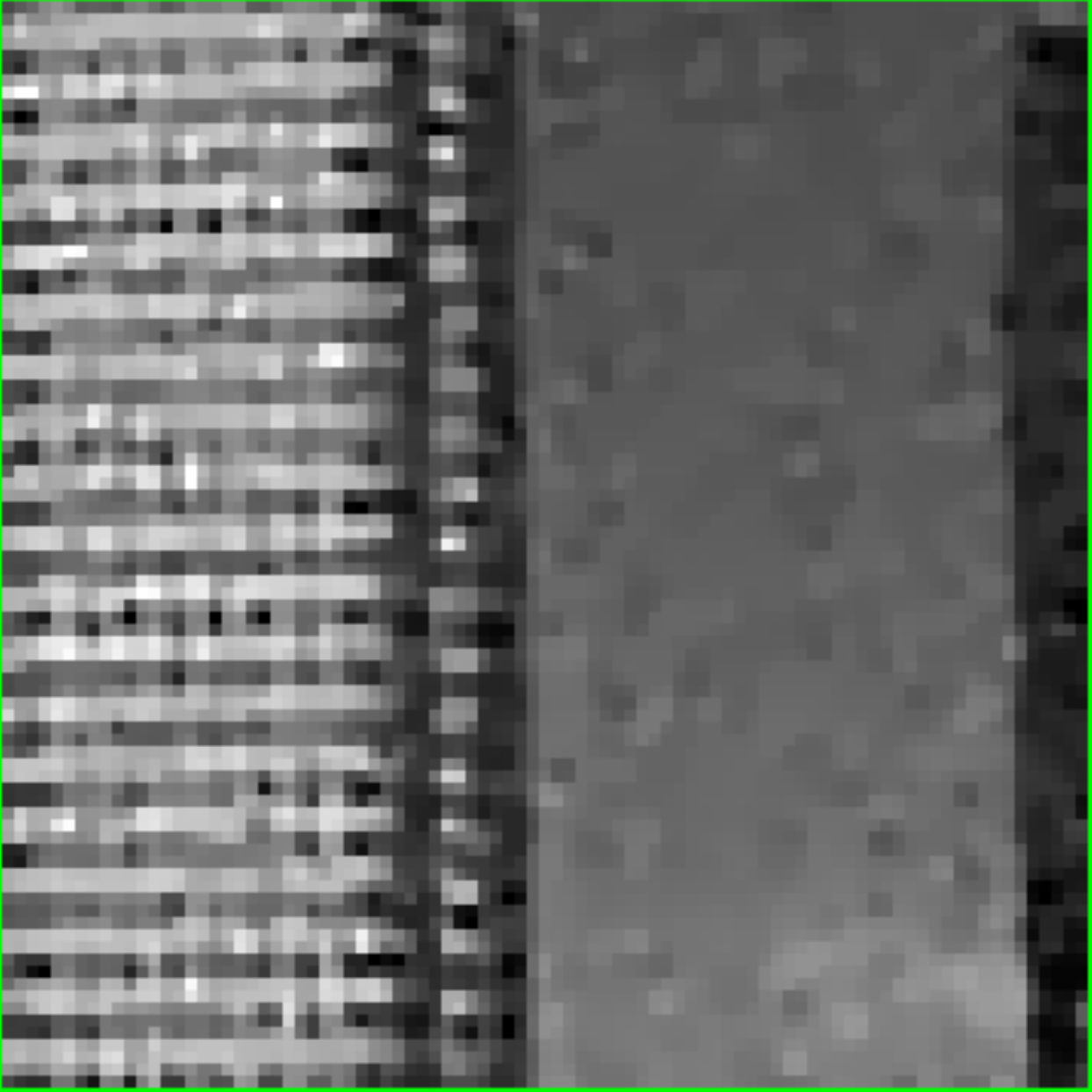}
		\caption{SATV (zoom)}
		\label{fig:SATVsky_zoom}
		\end{subfigure}
		\begin{subfigure}{0.24\textwidth}
		\centering
		\includegraphics[scale = 0.25]{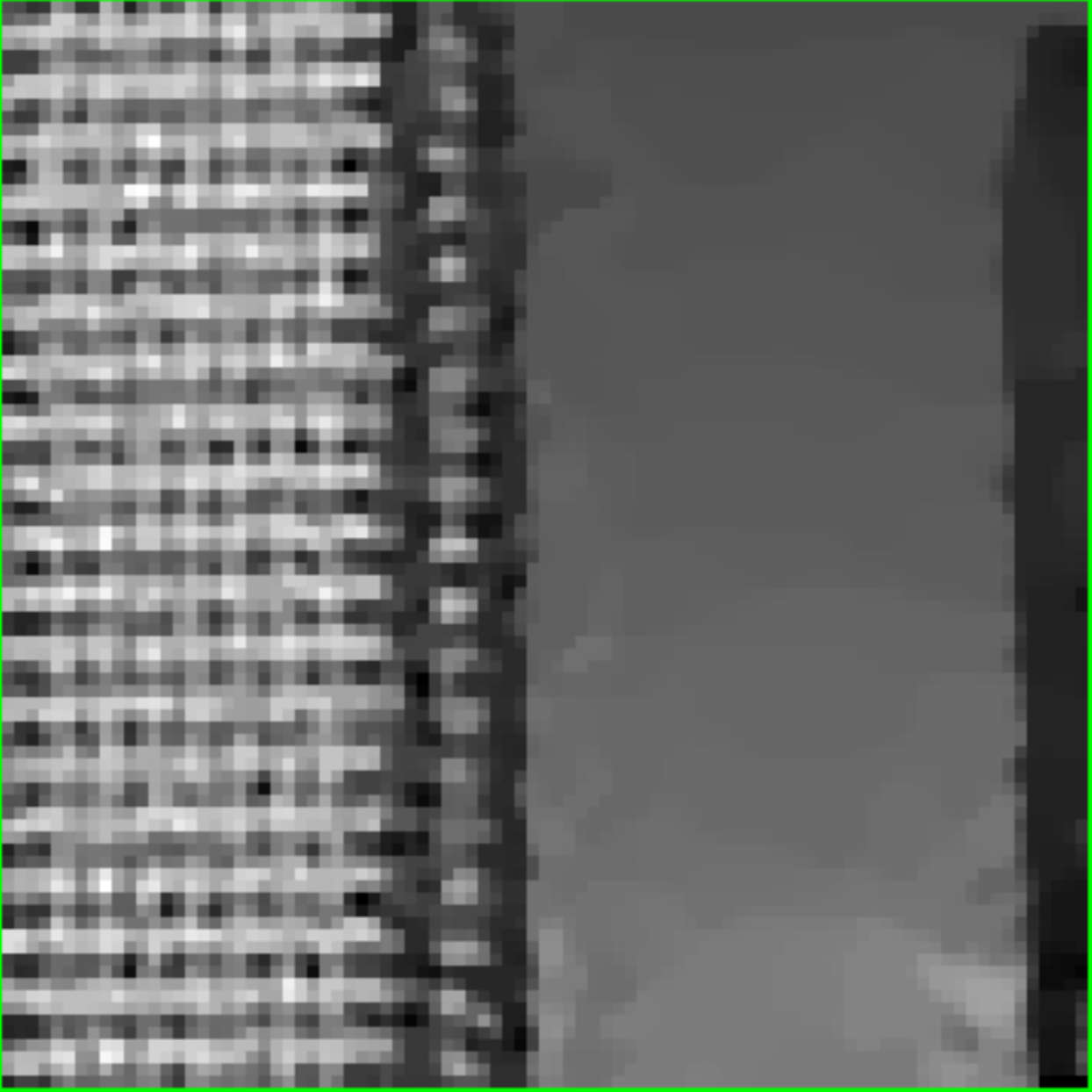}
		\caption{HWTV-L$_2$ (zoom)}
		\label{fig:WTVsky_zoom}
		\end{subfigure} 
		\begin{subfigure}{0.24\textwidth}
		\centering
		\includegraphics[scale = 0.25]{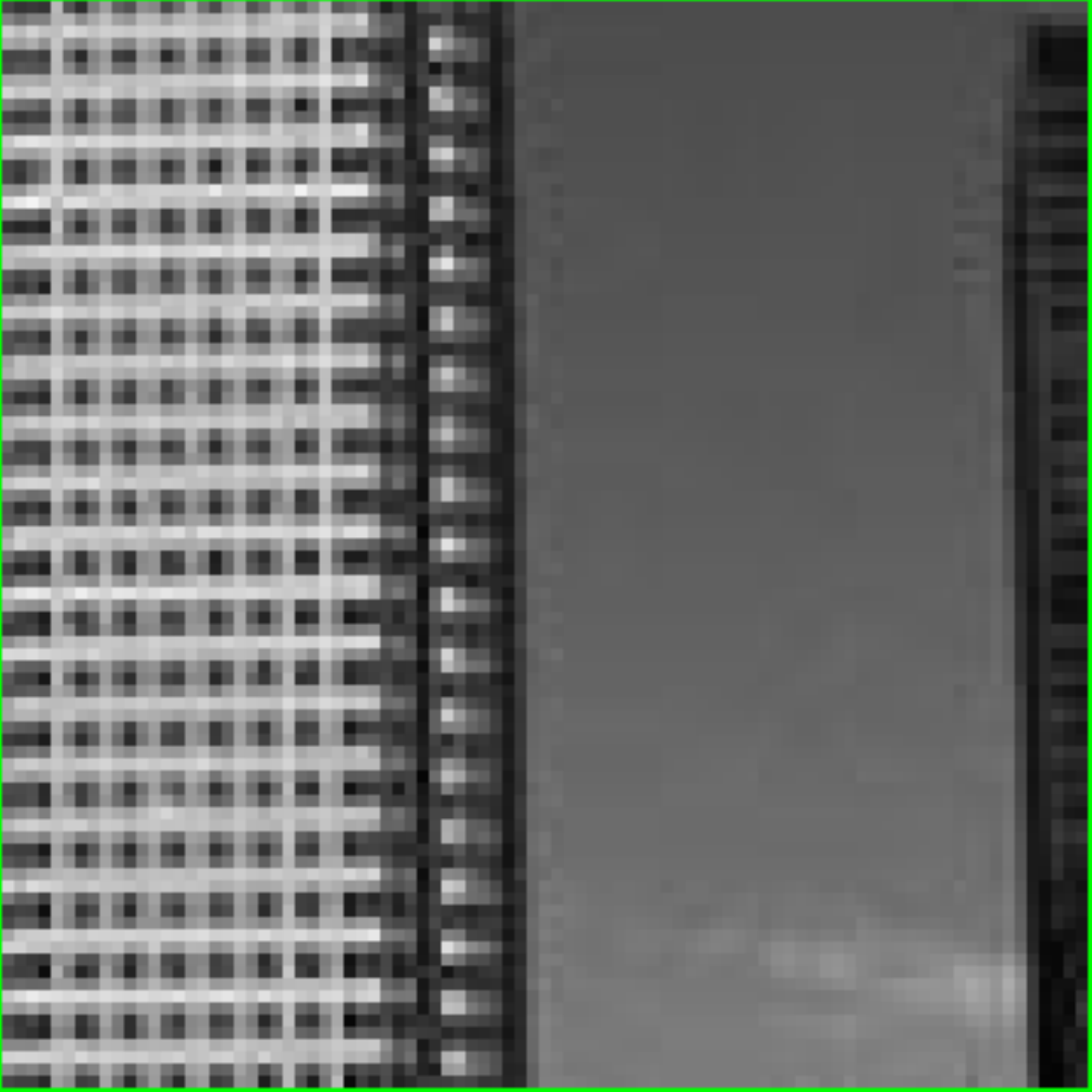}
		\caption{$u$ (zoom)}
        \label{fig:sky_zoom}
		\end{subfigure}
  \caption{\textbf{ISNR optimisation}. \emph{First row}:  Reconstruction of image in Fig.\ref{fig:skyscraper_or} by TV-L$_2$ ($\tau = 0.91$) \ref{fig:TVsky}, SATV ($\tau = 0.94$) \ref{fig:SATVsky}, HWTV-L$_2$ ($\tau = 0.94$, $r=14$) \ref{fig:WTVsky} and $\alpha$ parameter map obtained by the proposed parameter procedure \ref{fig:alpha}. \emph{Second row}: zoomed details.}
    	\label{fig:skyisnr}
\end{figure}

\begin{figure}[!t]
	\centering
		\begin{subfigure}{0.24\textwidth}
		\centering
		\includegraphics[scale = 0.34]{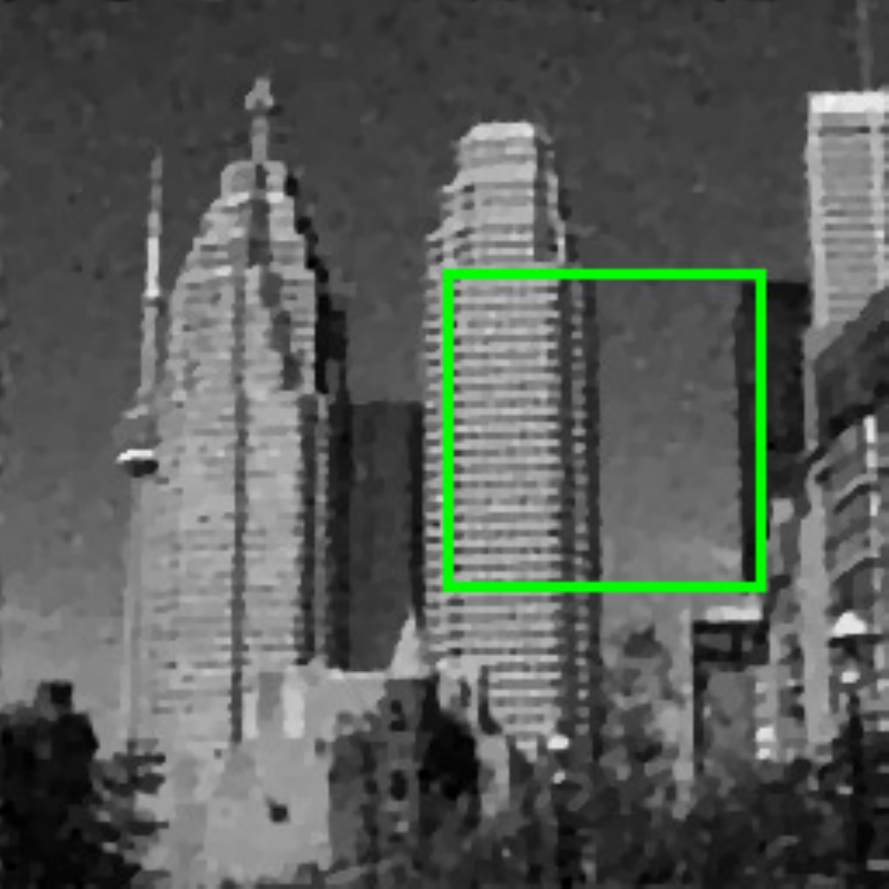}
		\caption{TV-L$_2$}
		\label{fig:TVskySSIM}
		\end{subfigure}
		\begin{subfigure}{0.24\textwidth}
		\centering
		\includegraphics[scale=0.34]{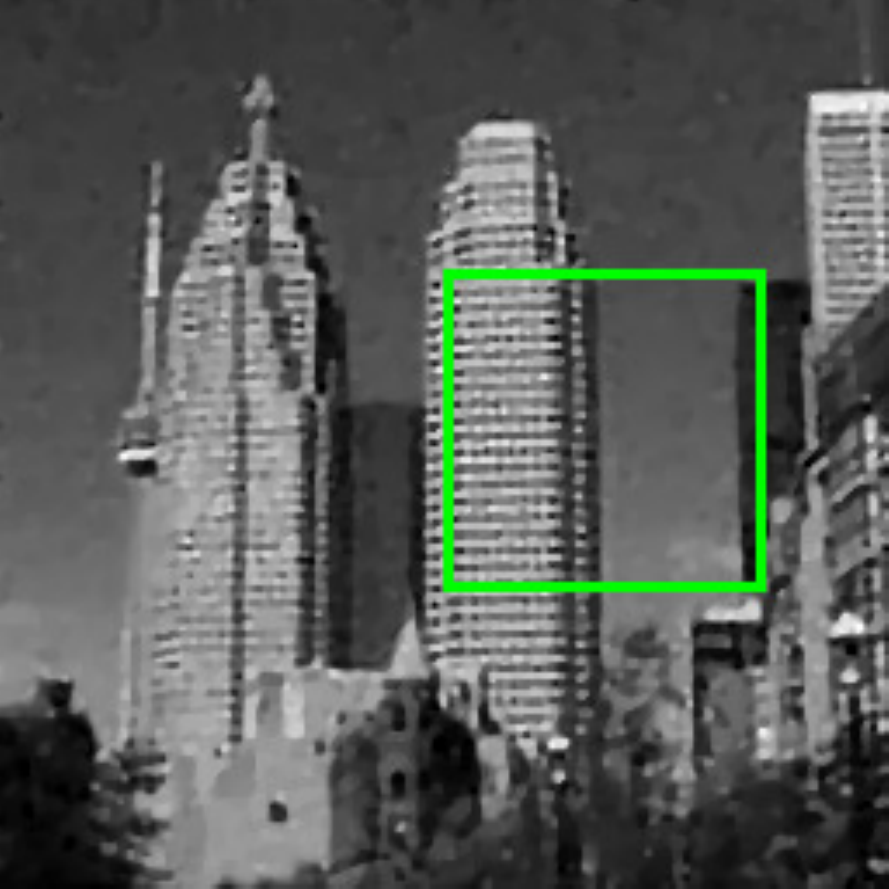}
		\caption{SATV}
		\label{fig:SATVskySSIM}
		\end{subfigure}
		\begin{subfigure}{0.24\textwidth}
		\centering
		\includegraphics[scale = 0.34]{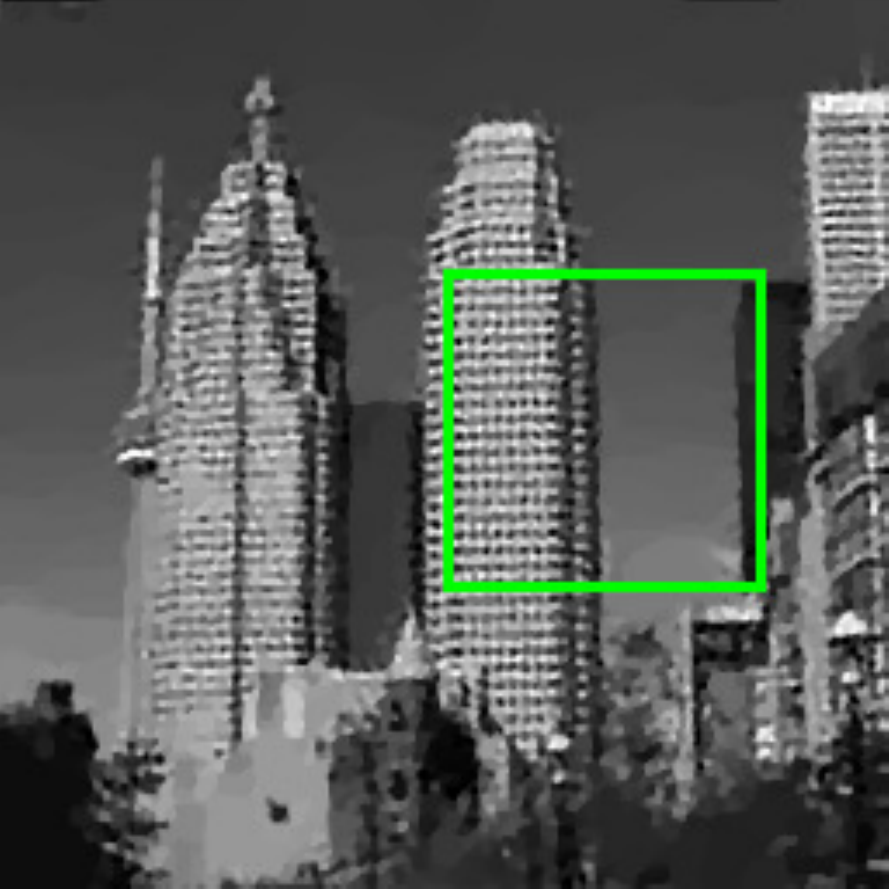}
		\caption{HWTV}
		\label{fig:WTVskySSIM}
		\end{subfigure} 
		\begin{subfigure}{0.24\textwidth}
		\centering
		\includegraphics[scale = 0.20]{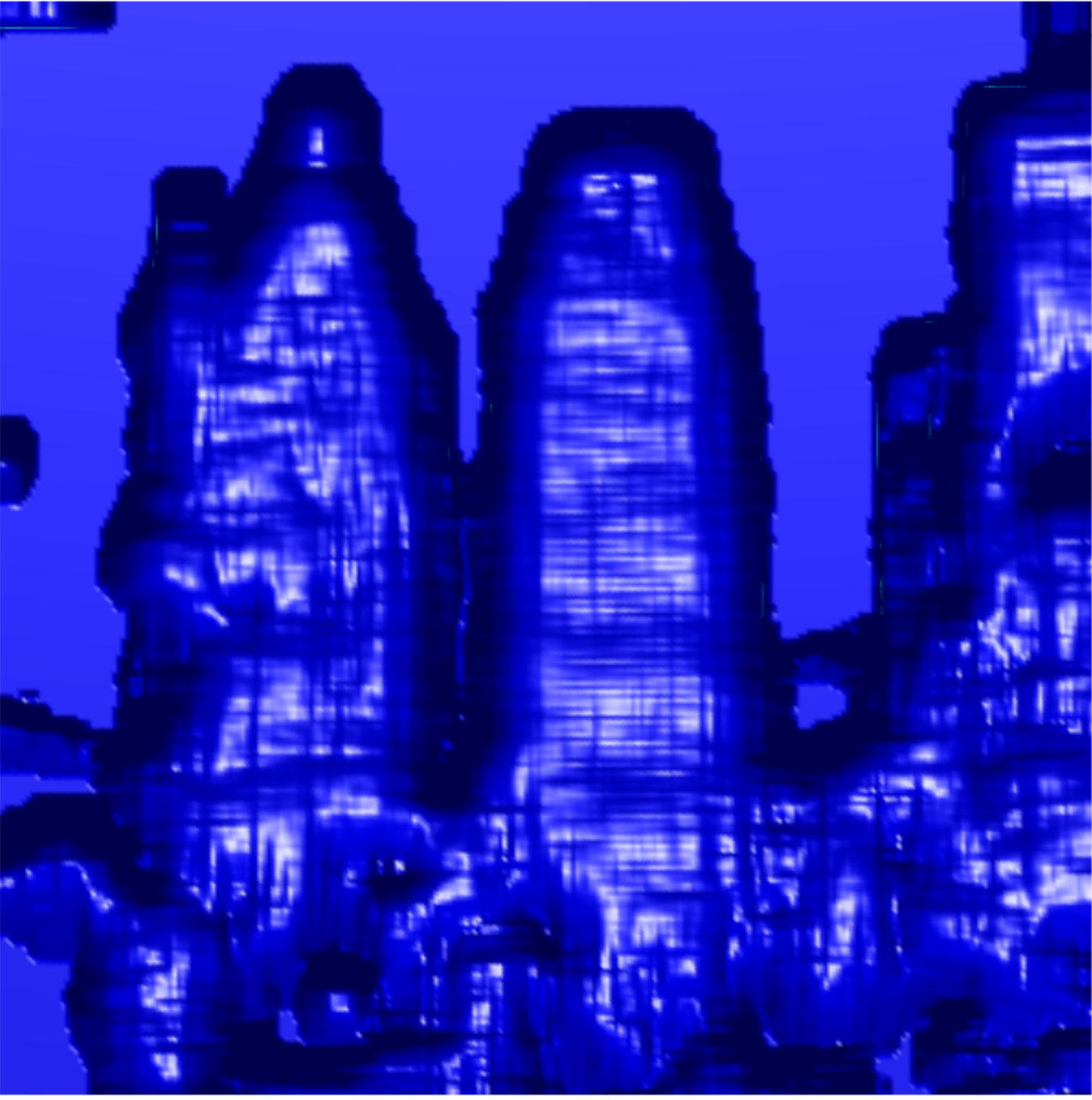}
		\caption{$\alpha$}
        \label{fig:alphaSSIM}
		\end{subfigure}\\
		\begin{subfigure}{0.24\textwidth}
		\centering
		\includegraphics[scale = 0.25]{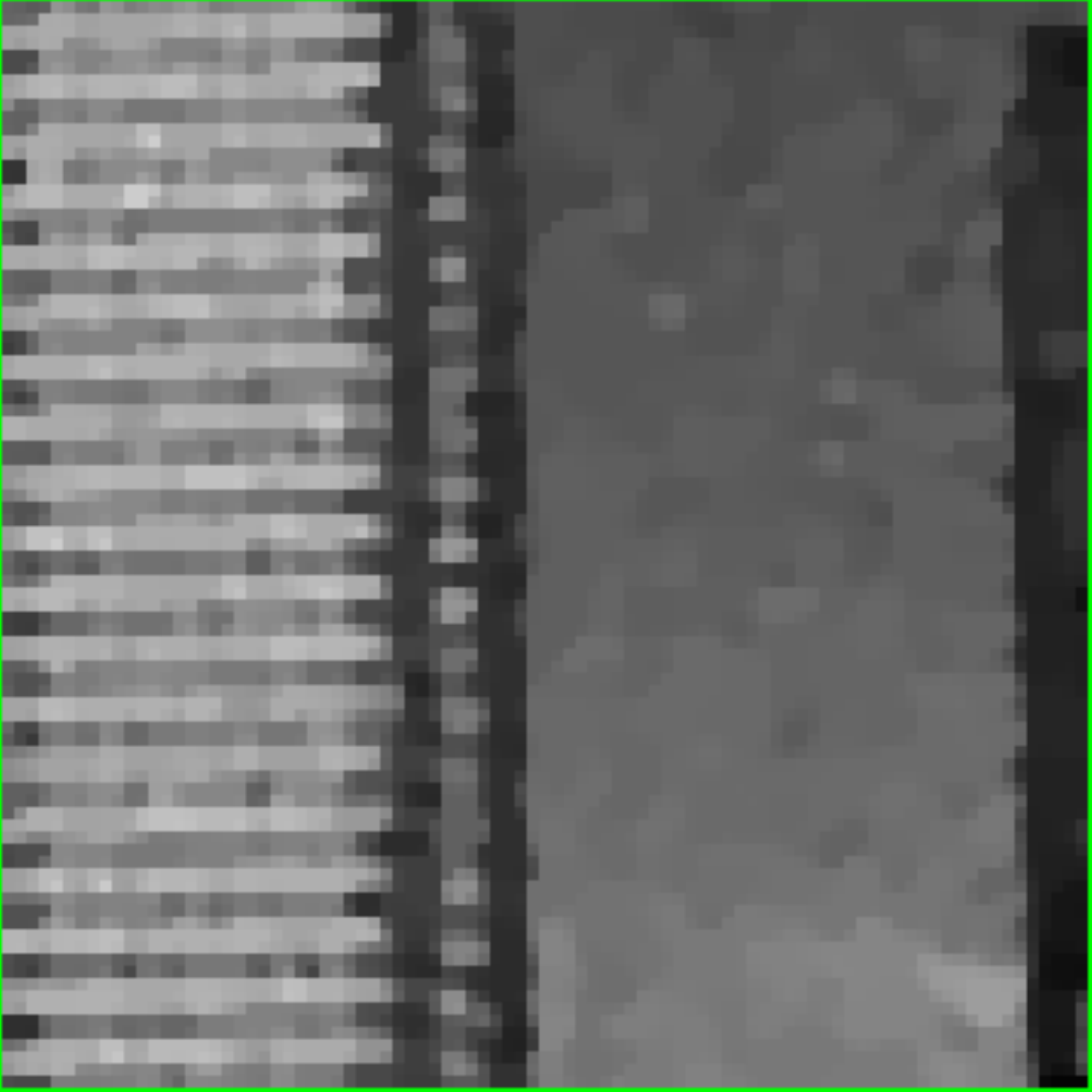}
		\caption{TV-L$_2$ (zoom)}
		\label{fig:TVsky_zoomSSIM}
		\end{subfigure}
		\begin{subfigure}{0.24\textwidth}
		\centering
				\includegraphics[scale=0.25]{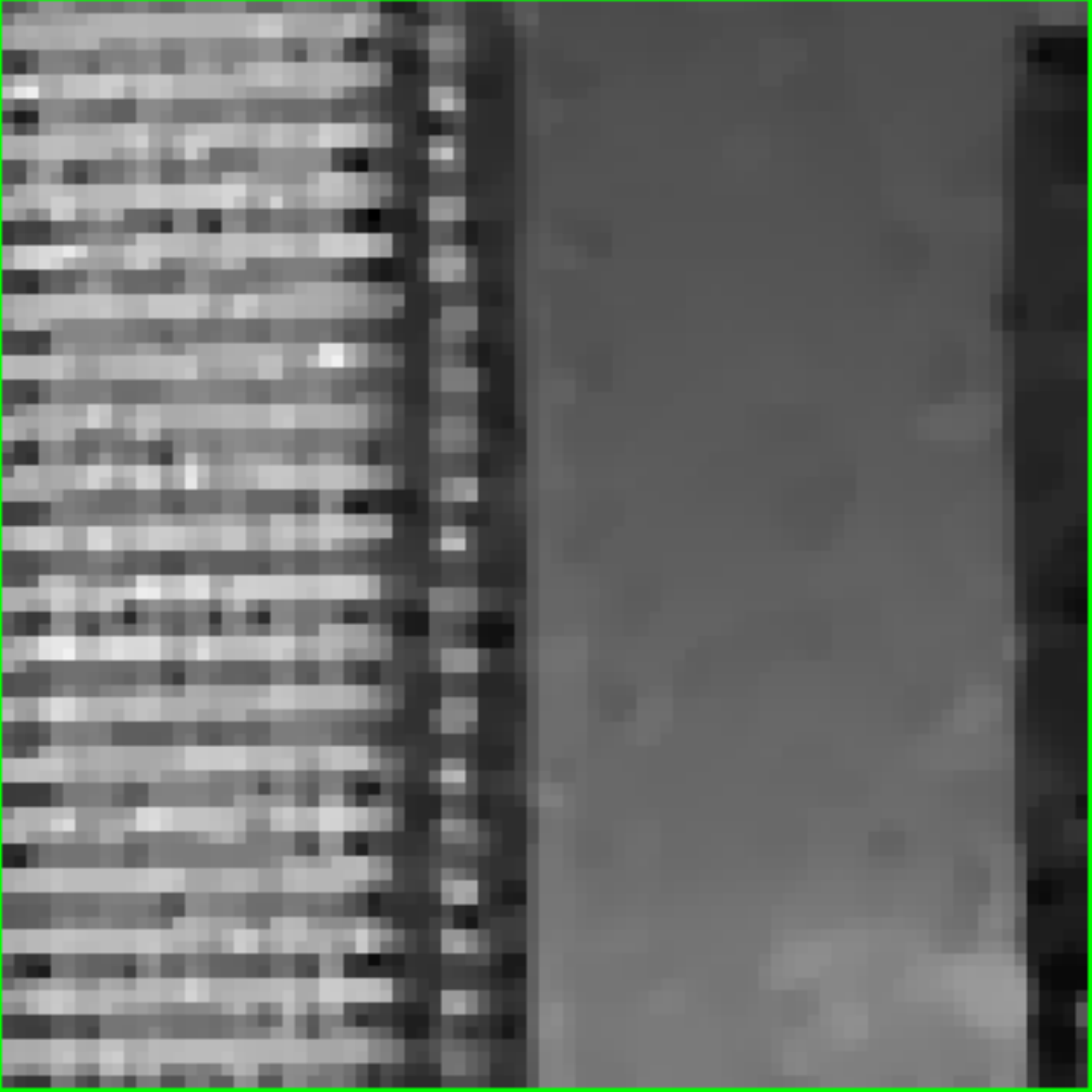}
		\caption{SATV (zoom)}
		\label{fig:SATVsky_zoomSSIM}
		\end{subfigure}
		\begin{subfigure}{0.24\textwidth}
		\centering
		\includegraphics[scale = 0.25]{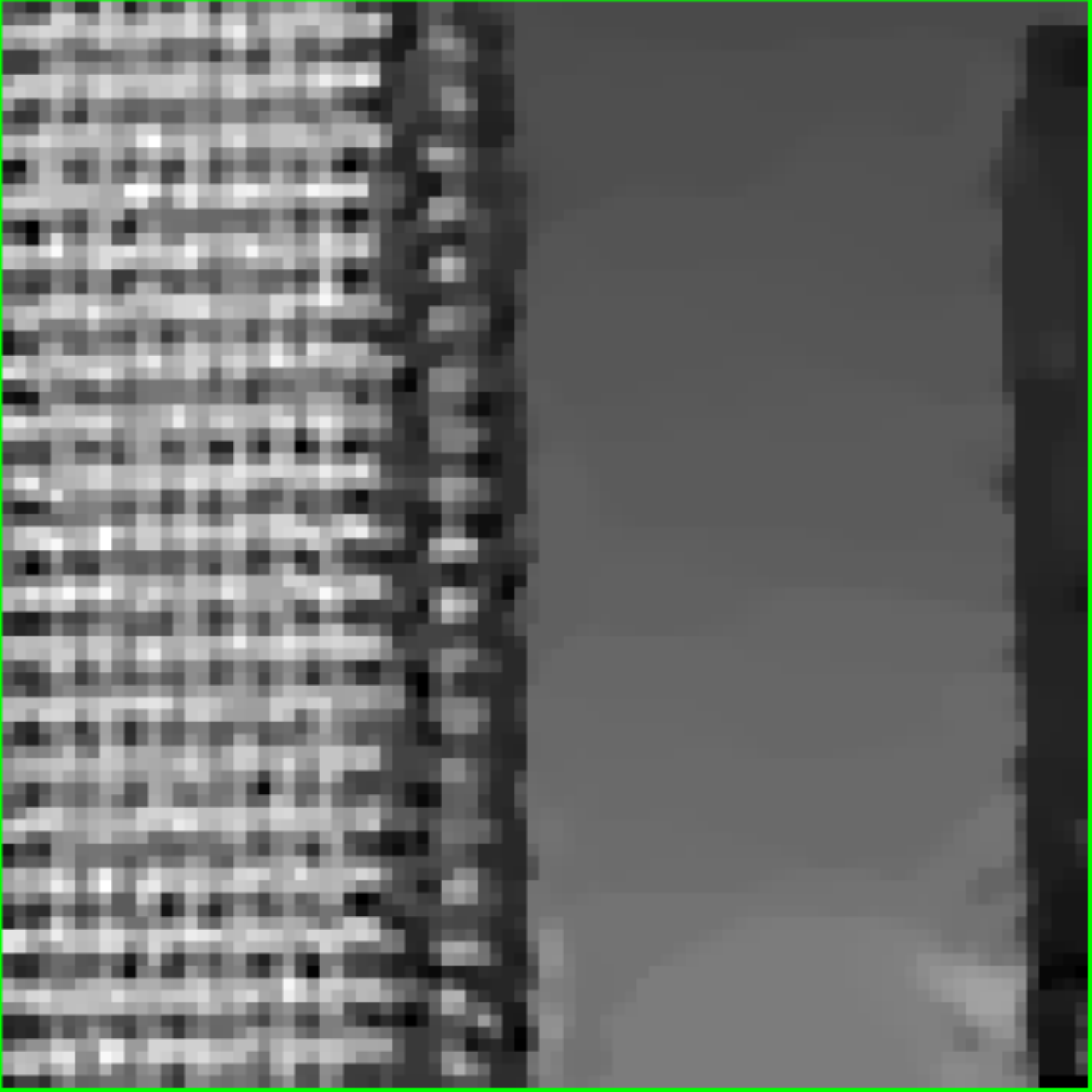}
		\caption{HWTV-L$_2$ (zoom)}
		\label{fig:WTVsky_zoomSSIM}
		\end{subfigure} 
		\begin{subfigure}{0.24\textwidth}
		\centering
		\includegraphics[scale = 0.25]{images/true_zoom.eps}
		\caption{$u$ (zoom)}
        \label{fig:sky_zoomSSIM}
		\end{subfigure}
  \caption{\textbf{SSIM optimisation}. \emph{First row}: Reconstruction of image in Fig.\ref{fig:skyscraper_or} by TV-L$_2$ ($\tau = 0.98$) \ref{fig:TVskySSIM}, SATV ($\tau = 0.95$) \ref{fig:SATVskySSIM}, HWTV-L$_2$ ($\tau = 0.93$, $r=6$) \ref{fig:WTVskySSIM} and $\alpha$ parameter map for WTV \ref{fig:alphaSSIM}. \emph{Second row}: zoomed details.}
    	\label{fig:skyssim}
\end{figure}

\paragraph{Image deblurring.} We consider the \texttt{skyscraper} test image ($256 \times 256$) corrupted by AWGN of two levels $\sigma=0.02,0.05$, and Gaussian blur of \texttt{band} = $5$ and \texttt{sigma} = $1$. The ground-truth image $u$ and the observed image $g$ for $\sigma=0.05$ are shown in Fig.\ref{fig:skyscraper_or}-\ref{fig:skyscraper_cor}, respectively. In this test we highlight the improvements obtained by our \eqref{eq:HWTV-L2} Algorithm \ref{alg:ADMM} in comparison to the standard \eqref{eq:TV-L2} model and the SATV method, both solved by means of ADMM for comparisons. As mentioned above, for the automatic adaption of the parameter $\mu$ along the iterations, a value for the parameter $\tau$ needs to be chosen. For the three models considered, we observed that the value of $\tau$ maximising the ISNR does not necessarily correspond to the value maximising the SSIM, see Fig.\ref{fig:ISNR_TVSATV}. For the TV-L$_2$ model the maximum SSIM is reached for $\tau\approx 1$, while the ISNR achieves its maximum when $\tau\approx 0.9$, the latter being the case in which texture is better preserved but noise is not completely removed. For SATV, the maximum ISNR and SSIM values are reached approximately for the same $\tau$. As remarked in \cite{satv}, the SATV method is robust w.r.t. the choice of the radius $r$ of the neighbourhoods used for the estimation. Thus, we set such parameter as the default value $r=5$ in our tests. We performed similar sensitivity tests for our HWTV-L$_2$ model for different $(\tau,r)$ values. Results are shown in Figs.\ref{fig:ISNR_SV}-\ref{fig:SSIM_SV}. 
For each method, we then selected the parameter(s) yielding the maximum ISNR/SSIM values and compared the results obtained. 
In Table \ref{tab:1} we report the achieved ISNR/SSIM values, whereas in Figs.\ref{fig:skyisnr}-\ref{fig:skyssim} we show the associated restored images for the case of AWGN with $\sigma = 0.05$ - see Fig.\ref{fig:skyscraper_or}(bottom). 
We observe that our HWTV-L$_2$ method results in higher quality reconstructions w.r.t. TV-L$_2$ and SATV. Visual inspection confirms the effectiveness of our approach in distinguishing between textured and homogeneous regions, see Figs.\ref{fig:WTVsky_zoom}-\ref{fig:WTVsky_zoomSSIM}.

\begin{figure}[!t]
	\centering
	\begin{subfigure}{0.24\textwidth}
	\centering
		\includegraphics[scale = 0.5]{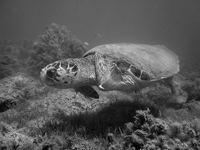}
	\caption{$u$}
	\label{fig:turtle_or}
	\end{subfigure}
	\begin{subfigure}{0.24\textwidth}
	\centering
		\includegraphics[scale =0.5]{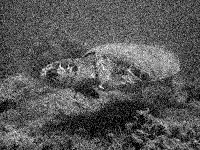}
		\caption{$g$}
		\label{fig:turtle_noisy}
	\end{subfigure}
	\begin{subfigure}{0.24\textwidth}	
	\centering
		\includegraphics[scale = 0.5]{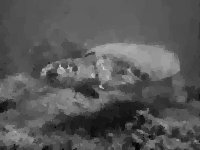}
		\caption{Bilevel WTV}
		\label{fig:turtle_bilevel}
	\end{subfigure}
	\begin{subfigure}{0.24\textwidth}
	\centering
		\includegraphics[scale = 0.5]{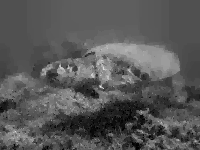}
		\caption{HWTV}
		\label{fig:turtle_ML}
    \end{subfigure}
\caption{\textbf{SSIM optimisation}. \ref{fig:turtle_or} Original image. \ref{fig:turtle_noisy} observed image corrupted by AWGN with $\sigma = 0.1$. \ref{fig:turtle_bilevel} WTV reconstruction obtained by bilevel optimisation of parameters $\alpha_i$ as in \cite{Hint1,Hint2} (SSIM = 0.7602). \ref{fig:turtle_ML} HWTV reconstruction ($\tau=0.86$, $r=40$, SSIM = \textbf{0.7708}).}
\label{fig:turtle}
\end{figure}

\paragraph{Image denoising.}
We now consider the test image \texttt{turtle}\footnote{Photo courtesy of K. Papafitsoros.} (150$\times$200) corrupted by AWGN of level $\sigma=0.1$ (see Figs.\ref{fig:turtle_or} -\ref{fig:turtle_noisy}) and focus on the quality and computational improvements of our HWTV-L$_2$ method in the case of anisotropic TV (i.e. $p=1$ in \eqref{TV:def}) in comparison to the alternative bilevel optimisation strategy used \cite{Hint1,Hint2} for estimating the space-variant parameters $\alpha_i$. 
After optimising the HWTV-L$_2$ method over $\tau$ as discussed above, the maximum achieved value is SSIM = 0.7708 (for $r=40$, $\tau=0.86$), in comparison to SSIM = 0.7602 obtained by using a bilevel optimisation strategy.
The reconstructions are shown in Fig.\ref{fig:turtle_bilevel}-\ref{fig:turtle_ML}.
We remark that, in addition to the obtained SSIM and visual improvements, our approach exhibits a very high computational efficiency, whereas bilevel codes are known to be computational expensive and hardly applicable to high-resolution images. For instance, in this experiments the proposed ADMM Algorithm \ref{alg:ADMM} for the HWTV-L$_2$ model required 
only 40 seconds on a standard laptop, compared to the 1429 seconds required by the bilevel algorithm \cite{Hint2}.
\begin{table}[!t]
\centering
	\begin{tabular}{c|ccc|cccc}
	    \hline\hline
        & &$\sigma=0.02$&&&&$\sigma=0.05$&\\
        \hline\noalign{\smallskip}
		&TV&SATV&HWTV&&TV&SATV&HWTV\\
		\hline
		ISNR&3.4701&3.6625&\textbf{4.3331}&&1.9433&2.0414&\textbf{2.5408}\\
		SSIM&0.8733&0.8966&\textbf{0.9007}&&0.7335&0.7797&\textbf{0.8099}\\
	\noalign{\smallskip}\hline\noalign{\smallskip}
	\end{tabular}
\caption{Maximum ISNR/SSIM values achieved by TV-L$_2$, SATV and HWTV-L$_2$ on the \texttt{skyscraper} test image in Fig.\ref{fig:skyscraper_or}(top) corrupted by AWGN of two different levels.}
	\label{tab:1}
\end{table}

\section{Conclusions}
We proposed an image restoration method based on an hybrid, locally-weighted TV-L$_2$ variational model where local regularisation parameters are combined with a global data fidelity weight.
Numerically, we solve the model by means of an ADMM-type algorithm combined with an effective and efficient automatic update of the local parameters via ML estimation and a global discrepancy constraint. Compared to standard as well as state-of-the-art competing models, the proposed approach outperforms in terms of standard image quality measures (ISNR, SSIM) as well as computational efficiency.

\section*{Acknowledgements}

LC acknowledges the support of the Fondation Math\'ematique Jacques Hadamard (FMJH). LC and MP are thankful to the organisers of the special trimester \emph{The mathematics of imaging} held at the IHP (Paris, France) where part of this research was carried out and to G. Peyr\'e for the financial support provided within the ERC project NORIA. Research of AL, MP and FS was supported by the ``National Group for Scientific Computation
(GNCS-INDAM)'' and by the ex60 project ``Funds for selected research topics''. The authors are deeply grateful to K. Papafitsoros for the computation of the bilevel WTV reconstruction used for comparison in Fig. \ref{fig:turtle}.

\end{document}